\documentclass[proceedings%
]{aofa}

\usepackage[utf8]{inputenc}
\usepackage{subfigure}

{\newtheorem{definition}{Definition}
\newtheorem{lemma}{Lemma}
\newtheorem{theorem}{Theorem}
\newtheorem{proposition}{Proposition}
\newtheorem{corollary}{Corollary}}

{
\newtheorem{example}{Example}
\newtheorem{remark}{Remark}
}
\def\BeginExample{\begin{example}}
\def\EndExample{\QED\end{example}}

\def\BeginRemark{\begin{remark}}
\def\EndRemark{\QED\end{remark}}

\def\BeginProof{\noindent{\it Proof\hskip1pt}: }
\def\BeginProofOf#1{\noindent{\it Proof of #1\hskip1pt}: }
\def\EndProof{\QED\betweenskip}
\def\betweenskip{\vskip10pt}

\def\boxitat#1#2{\vbox             
     {\hrule\hbox{\vrule\kern#1%
     \vbox{\kern #1\hbox{#2}\kern#1}\kern#1\vrule}\hrule}}
\def\enclose#1{\boxitat{0pt}{#1}}

\def\qedspace{\null}
\def\qedblack{\qedspace                    
  \lower.6pt\hbox{\vrule height7pt width 5pt}}
\def\qedwhite{\qedspace                    
  \lower.6pt\enclose{%
        \hbox{\vrule height7pt width 0pt\hskip6pt}}}
\def\QED{\qedwhite}

\def\beq{\begin{equation}}
\def\eeq{\end{equation}}

\usepackage{amssymb} 

\def\naturals{{\mathbb N}}

\def\eps{{\varepsilon}}

\def\sS{{\cal S}}
\def\sP{{\cal P}}
\def\sA{{\cal A}}
\def\sW{{\cal W}}
\def\sR{{\cal R}}

\def\sF{{\cal F}}

\def\sL{{\cal L}}

\def\kbar{{\bar k}}
\def\cvec#1{{\mathbf{#1}}}

\def\part#1{#1}
\def\cvec#1{{\mathbf{#1}}}
\def\term{{\cal F}}
\def\finish{finish }
\def\free{{\cal L}}
\def\span{{m}}

\def\len{{\rm len}}
\def\SEQ{{\rm SEQ}}

\def\Ex{{\mathbb E}}

\def\pr{{\mathbb P}}

\author[E.A. Bender and Z. Gao]{Edward A. Bender\addressmark{1}
\and Zhicheng Gao\addressmark{2}\thanks{Research supported by NSERC}}

\title[Restricted Sequential Structures and Runs]
{Locally Restricted Sequential Structures and Runs of a Subcomposition in Integer Compositions}

\address{\addressmark{1}Department of Mathematics, University of California, San Diego, La Jolla, CA 92093-0112\\
\addressmark{2}School of Mathematics and Statistics, Carleton University, Ottawa, Ontario K1S5B6
}

\begin{document}
\maketitle


\begin{abstract}
We study part sizes of supercritical locally restricted sequential
structures. This extends previous results about locally restricted
integer compositions and part sizes in smooth supercritical
compositional structures. Applications are given for runs of
subcompositions. The problems are formulated as enumerating
directed walks in sized infinite digraphs and the proofs depend heavily on
earlier results by Bender and Canfield about infinite transfer matrices.
\end{abstract}
\section{Introduction}

In \cite{BG14}  part sizes of compositional structures  were studied. It was shown that if
the composition is smooth supercritical then the numbers of parts of large sizes are asymptotically Poisson and
an asymptotic expression was obtained for the expected value of the maximum part size which is accurate up to $o(1)$. In \cite{BCG}, part sizes of locally restricted nearly free integer compositions were studied. We will extend some of the major results of \cite{BCG,BG14} to locally restricted supercritical sequential structures.
Runs of a single letter/part in words/compositions have been studied extensively; see e.g. \cite{BG14,Gafni,Wilf}.   In locally restricted structures, such as Carlitz compositions and Smirnov words, runs of a single part may not be allowed, and hence it is natural to consider runs of a substructure. Our main results on large part size distributions will be applied to maximum run length of a subcomposition in several classes of locally restricted compositions. Our approach follows that in \cite{BC09,BCG} using directed walks in sized infinite digraphs and properties of the corresponding infinite transfer matrices.
Our proofs rely heavily on
results in \cite{BC09,BCG} about infinite transfer matrices and large part size distributions.

\section{Definitions}

\begin{definition}[Sequential structures]
Let $\sP$ be a class of combinatorial structures, called {\em parts}.
\begin{itemize}
\item
Each part $\part p$ has a positive integer {\em size},
denoted $|\part p|$.
\item
We use $\sP_n$ to denote the set of parts of size $n$ and
assume $P_n:=|\sP_n|<\infty$ for each $n\in\naturals$.
\item
For each integer $k\ge0$, we denote the class consisting of sequences
of $k$ parts by
$$
\SEQ_k(\sP):=\{\part p_1\part p_2\ldots\part p_k: \part p_j\in\sP\}.
$$
Thus $\SEQ_0(\sP)$ contains only the empty sequence which has size and length~0.
{\bf Throughout this paper} $\eps$ denotes the empty sequence.
\item
Let $\SEQ_{<k}(\sP)$ denote the class of  sequences of at
most  $k-1$ parts.\\
Let $\SEQ(\sP):=\cup_{k\ge 0}\SEQ_k(\sP)$,
the set of all sequences.
\item
If $\cvec a=\part p_1\part p_2\ldots\part p_k \in \SEQ(\sP)$,
then the {\em length} of $\cvec a$ is $k$ and we
write $\len(\cvec a)=k$.
The {\em size} of $\cvec a$ is
$$
|\cvec a|=|\part p_1|+|\part p_2|+\cdots+|\part p_k|
$$
and the {\em distance} from $\part p_i$ to $\part p_j$ is $|i-j|$.
\end{itemize}
\end{definition}

\BeginExample[Generalized compositions]\label{ex:gen-comp}
When studying $\SEQ(\sP)$, only the values of $P_n$ are important.
Thus we could think of $\sP$ as parts in a generalized composition
where the part $n$ comes in $P_n$ ``colors'', or whatever you choose
to call them.
If $P_n=0$, there are no parts of size~$n$.
Thus, ordinary compositions correspond to $P_n=1$ for all $n$ and
words on a $k$-letter alphabet correspond to $P_1=k$ and $P_n=0$
for $n>1$.
When parts of size $n$ occur in $n$ colors, $P_n=n$ and we have $n$-colored
compositions, which were studied in~\cite{Ag00,NA06,NA08}. We may also consider colored
compositions where a part of size $n$ corresponds to a multiset of $n$ colored balls with
$N$ colors available. Here we have $P_n={n+N-1\choose N-1}$.
\EndExample

Let $\sA\subseteq \SEQ(\sP)$. It is called locally restricted if the parts of a structure in $\sA$ within a fixed distance satisfy certain restrictions.
Locally restricted integer compositions were studied in \cite{BC09}, where local restrictions are defined in terms of local restriction functions.
The function was then used to construct a digraph.
In this paper, we define local restrictions directly in terms of the digraph.
Readers wishing to see the connection between local restriction functions and the digraph should consult~\cite{BC09}.

\begin{definition}
[Locally restricted structures]
\label{def:LRS}
Let $\span\in \naturals$, $\sS,\term\subseteq\SEQ_{<\span}(\sP)$ and $\sR\subseteq\SEQ_{\span}(\sP)$. The integer $m$
is called the {\em span} of the locally restricted class of structures associated with the digraph $D$ which has vertex set $V(D)=\sS\cup\sR\cup\term$.
If $p\in\sS\cap\term$, we allow two copies as separate vertices in the digraph. For convenience, we introduce two copies of the empty sequence, denoted by $\eps_s$ and $\eps_f$,  such that $\eps_s\in \sS$ and $\eps_f\in \term$.
Suppose $D$ satisfies the following conditions.
\begin{itemize}
\item[(a)] There is an arc from $\eps_s$ to every other vertex in $\sS$, and at least one arc from $\sS$ to $\sR$.
\item[(b)] There is an arc to $\eps_f$ from every other vertex in $\term$, and at least one arc from $\sR$ to $\term$.
\item[(c)] The sub-digraph $D_{\sR}$ of $D$ induced by $\sR$ is strongly connected, $|\sR|\ge 2$, and the sub-digraph of $D$ induced by $\sS\cup \sF$ contains no directed cycle.
\end{itemize}

The vertices of $\sR$ will be called the {\em recurrent vertices},
and the vertices of $\sS$ and $\term$ will be called the
{\em  start vertices} and {\em  \finish vertices}, respectively.
If a part $\part p\in \sP$ appears in some recurrent vertex, then
$\part p$ is called a {\em recurrent part}.

Let $\sW$ denote the set of all directed walks in $D$ from
$\eps_s$ to $\eps_f$.
We use $\SEQ(\sP;D)$ to denote the class of all structures of the concatenation form
$\cvec v_1\cvec v_2\cdots\cvec v_j$, where $\eps_s\cvec v_1\cvec v_2\cdots\cvec v_j\eps_f\in \sW$.
If each element in $\SEQ(\sP;D)$ arises from {\bf only one} walk in $\sW$,
we say that $\SEQ(\sP;D)$ is a {\em locally restricted class} associated
with $D$.
\end{definition}
\noindent
The ``only one'' condition is required so that an ogf built from $D$
will count each element in $\SEQ(\sP;D)$ just once.
The span is actually associated with the digraph.
One can easily construct a digraph $D'$ with span $k\span$ for any $k\in\naturals$
such that $\SEQ(\sP;D')=\SEQ(\sP;D)$.

\begin{definition}[Regular class]\label{def:regular}
A class $\SEQ(\sP;D)$ of locally restricted structures will be called {\em regular} if it satisfies the following conditions.
\begin{itemize}
\item
The gcd of the lengths of all directed cycles in $D_{\sR}$ is equal to 1.
\item
There is a positive integer $k$ and vertices $\cvec v_0, \cvec v_k\in \sR$ such that
$\gcd\{m-n:m,n\in S\}=1$, where
$S=\{n:n=|\cvec v_0|+\cdots+ |\cvec v_k| \hbox{ for some directed walk $\cvec v_0\cdots\cvec v_k$ of length $k$ in $D_{\sR}$} \}$.
\end{itemize}
\end{definition}

\BeginExample[Pattern avoidance]
Let $B$ be a finite set of structures in $\SEQ(\sP)$, and $\sA$ be the class of structures in $\SEQ(\sP)$ which don't contain any structure in $B$. We say that the structures in $\sA$ {\em avoid} the structures in $B$ (or simply, avoid $B$). We may construct the class $\sA$ using the following digraph. Let $m+1$ be the maximum length of the structures in $B$. Let $\sS=\{\eps_s\}$, $\sF\subset \SEQ_{<m}(\sP)$ and $\sR\subset \SEQ_{m}(\sP)$ be consisting of structures which avoid $B$. There is an arc from vertex $\cvec a$ to a vertex $\cvec b$ if and only if $\cvec a\cvec b$ avoids $B$. There is also an arc from $\eps_s$ to every vertex in $\sR\cup\sF$. Then $\sA=\SEQ(\sP;D)$. We note that the second gcd condition in Definition~3 is not satisfied here because each recurrent vertex has size $m$.
\EndExample

\noindent {\bf Remark:} Words over a finite alphabet which avoid certain patterns have been studied extensively (see, e.g., \cite{HM09}). Here the size of a word is the length of the word. Since $\sP$ is finite here, one may use the simpler transfer matrix $T(z)$ such that $T_{i,j}(z)=z^{|r_j|}$ for recurrent vertex $\cvec r_j$ (See Definition~6 on page~8 for the transfer matrix). Consequently the ogf $F_R(z)=\cvec s^t(I-T(z))^{-1}\cvec f$ is a rational function of $z$. We note that in this case, each recurrent vertex has size $m$ and so $T(z)$ is a function of $z^m$, and hence we may apply \cite[Theorem 1]{BC09} to $T$ with $x=z^m$.

\begin{definition}[Generating functions and supercritical structures]
\label{def:super}
Let $\sA\subseteq\SEQ(\sP)$.
\begin{itemize}
\item
Define the ordinary generating function (ogf)
$P(z) := \sum_{p\in\sP}z^{|p|} =  \sum_{n\ge 1}P_nz^n$.
\item
{\bf Throughout this paper} $\rho$ is the radius of convergence of $P(z)$.
\item
If the radius of convergence of the ogf $A(z):=\sum_{\cvec a\in\sA}z^{|\cvec a|}$
is less than $\rho$, we call $\sA$ {\em supercritical}.
\end{itemize}
\end{definition}

\BeginExample[Alternating compositions]\label{ex:alternating}
Alternating, or up-down, compositions are compositions in which a part is
alternately greater and less than the preceding part.
 We set $\span=2$, $\sP=\naturals$,
$\sS=\SEQ_{<2}(\sP)$, $\term=\SEQ_{<2}(\sP)\backslash\{1\}$, and
$\sR=\{i,j\mid i>j\}$.
The empty sequence in $\sS$ (resp.\ $\term$) connects to every element in $\sR$.
There is an arc from $i\in\sS$ to $j\in\term$ whenever $i<j$ or $i$ is the empty sequence.
The arcs having at least one end in $\sR$ should be fairly easy to see.
Every alternating composition, including the empty one, is associated with a
unique directed walk from $\eps_s$ to $\eps_f$.  It is known \cite{BC09} that the class of alternating compositions is regular and supercritical.
\EndExample

\begin{definition}[Asymptotically free set]
\label{def:free}
Let $\sA=\SEQ(\sP;D)$ be a class of locally restricted structures with $|\sP|=\infty$
and let $\free$ be an infinite set of recurrent parts.
If $\cvec a=\part p_1\ldots\part p_k$ and $t=\min(k,\span-1)$, define $s(\cvec a):=\part p_1\ldots\part p_t$
and $f(\cvec a):=\part p_{k-t+1}\ldots\part p_k$.
Suppose there is a function $g:\SEQ_{<m}(\sP)\times\SEQ_{<m}(\sP)\to\naturals$ such that the following holds.
If $\cvec a l\cvec z\in \SEQ(\sP;D)$ with $l\in\sL$ and
$| l|\ge g(f(\cvec a),s(\cvec z))$, then $\cvec a l'\cvec z\in\SEQ(\sP;D)$ whenever $ l'\in\sL$ and $| l'|\ge g(f(\cvec a),s(\cvec z))$.
Then we call $\sL$ {\em asymptotically free}.
\end{definition}
What this says is that if $ l$ is large enough as determined by parts closer than distance $m$,
then it can be replaced by any large enough part in $\sL$.

\BeginExample[$k$-Carlitz compositions]\label{ex:m-Carlitz}
A $k$-Carlitz composition is a composition in which each part is different from each of the preceding $m$ parts,
i.e.\ $p_i\ne p_j$ whenever $|i-j|\le m$.
(Carlitz compositions, for which adjacent parts differ, is the case $k=1$.
They have been studied extensively.
See for example~\cite{LouPro02}.)
They are regular supercritical structures and a digraph $D$ can be constructed with any $\span\ge k$.
Let $\sP=\naturals$ and $\sS=\{\eps_s\}$.
Let $\sR$ (respectively $\term$) be all compositions in $\SEQ_\span(\sP)$ (respectively $\SEQ_{<\span}(\sP))$
which are $k$-Carlitz.
There is an arc from $\eps_s$ to all vertices in $\sR\cup \term$.
There is an arc from $\cvec a\in\sR$ to $\cvec b\in\sR\cup\term$
whenever $\cvec a\cvec b$ is an $k$-Carlitz composition.
Since a part that is larger than all parts within distance $k$ can be replaced by any larger part,
the set $\naturals$ is asymptotically free in $k$-Carlitz compositions.
\EndExample

\section{Main Results}

Let $r$ be the radius of convergence of the ogf for $\SEQ(\sP;D)$.
In the following, all logarithms will be to the base $1/r$.
Our main results are the following.

\begin{theorem}\label{thm:main}
Let $\sA=\SEQ(\sP;D)$ be a regular supercritical class of locally restricted structures,
and $\free=\{l_1,l_2,\ldots\}$ be an asymptotically  free set. Assume $|l_1|<|l_2|<\cdots$.
Select a structure $\cvec a$ uniformly at random from $\sA_n$. Let $\zeta_k(n)$ be the number of occurrences of $l_k$ in
$\cvec a$.
The following are true.
\begin{itemize}
\item[(a)] $\displaystyle |\sA_n|=Ar^{-n}\left(1+O\left(e^{-\delta n}\right)\right)$ for some positive constants $r$, $A$ and $\delta$.
\item[(b)] The distribution of $\zeta_k(n)$ is asymptotically normal with mean and
variance asymptotically proportional to $n$.
\item[(c)] The limit
\begin{equation}\label{eq:Ex/Ex}
v_k = \lim_{n\rightarrow\infty}
\frac{\Ex(\zeta_k(n))}{n}
\end{equation}
exists, and $v_k \sim Cr^{|l_k|}$ as $k\to\infty$, for some positive
constant $C$.
\item[(d)]
Suppose there is a function $\omega_1(n)\to \infty$  such that $\{|l_k|:k\ge 1\}\cap [\log n-\omega_1(n), n]$ is not empty for all sufficiently large $n$.
Then there is a function
$\omega_2(n)\rightarrow\infty$ such that the random variables
$\{\zeta_k(n):\log n -\omega_2(n) \le |l_k|\le n\}$ are asymptotically
independent Poisson random variables with means $\mu_k =Cnr^{|l_k|}$.
\end{itemize}
\end{theorem}

The distribution of large part sizes in general supercritical compositional structures has been studied extensively. Some latest results can be found in \cite{BG14}. We will convert runs of a given subcomposition into free parts in a related class of locally restricted structures and apply Theorem~1(d) to derive the following.
\begin{theorem}\label{thm:Run}
Let $\cvec c$ be a given composition such that $\cvec c$ cannot be written in the form $\cvec c=\cvec x\cvec y\cvec x$. Let $\sA=\SEQ(\naturals;D)$ be a regular supercritical class of locally restricted compositions with span $m=\len(\cvec c)$. Assume that $\cvec c$ is a recurrent vertex and the digraph $D$ contains the arc from  $\cvec c$ to itself. Let $C,r$, $\cvec a$ and $\zeta_k$ be as defined in Theorem~\ref{thm:main},
$R_n$ be the maximum run length of $\cvec c$ in $\cvec a$, and $g_n(k)$ be the probability that $\cvec a$ contains exactly $k$ runs of $\cvec c$ of length $R_n$. Let  $\gamma\doteq 0.577216$ be Euler's constant and define
\begin{equation}\label{eq:oscillation}
P_k(x) ~= ~\frac{\log e}{|\cvec c|}\sum_{\ell\ne 0}
\Gamma\left(k+\frac{2i\pi \ell\log e}{|\cvec c|}\right)\exp\left(\frac{-2i\ell\pi\log x}{|\cvec c|}\right).
\end{equation}
Then
\begin{itemize}
\item[(a)]
$\displaystyle \pr(R_n<k)\sim \exp\left(-\frac{Cn}{1-r^{|\cvec c|}}r^{k|\cvec c|}\right)$.
\item[(b)]
$\displaystyle \Ex(R_n)=\frac{1}{|\cvec c|}\log\frac{Cn}{1-r^{|\cvec c|}}+\frac{\gamma\log e}{|\cvec c|}-\frac{1}{2}-P_0\left(\frac{Cn}{1-r^{|\cvec c|}} \right)+o(1)$.
\item[(c)]
$\displaystyle g_n(k)=\frac{(1-r^{|\cvec c|})^k}{k!}P_k\left(\frac{Cn}{1-r^{|\cvec c|}}\right)+\frac{(1-r^{|\cvec c|})^k\log
e}{k|\cvec c|}+o(1)$.
\end{itemize}
\end{theorem}

\medskip
\begin{corollary}
Theorem~2 holds for the following classes of compositions.
\begin{itemize}
\item[(a)] All compositions for any given composition $\cvec c$, where $r=1/2$ and $C=\frac{1}{2}(1-2^{-|\cvec c|})^2$. In particular, $C=1/8$ when $\cvec c=1$, which gives Gafni's result \cite{Gafni}.
\item[(b)] Carlitz compositions for any given Carlitz composition $\cvec c$, where $r\doteq 0.571350$ is the smallest positive number satisfying
    $\displaystyle \sum_{j\ge 1}\frac{r^j}{1+r^j}=1$. In particular, when $\cvec c=ab$ with $a\ne b$, we have
     $$ C=\frac{\left(1-r^{a+b}\right)^2}{(1+r^a)(1+r^b)}\frac{1}{\sum_{j\ge 1}\frac{jr^j}{(1+r^j)^2}}.$$
\item[(c)] $k$-Carlitz compositions for any given $k$-Carlitz composition $\cvec c$.
\item[(d)] Alternating compositions for  any given alternating composition $\cvec c$, where $r\doteq 0.6363$.
\item[(e)] $n$-color compositions (defined in Example~1) for any given $n$-color composition $\cvec c$, where $r=\frac{3-\sqrt{5}}{2}$.
\item[(f)] Colored compositions with $P_n={n+N-1\choose N-1}$ (defined in Example~1) for any given colored composition $\cvec c$, where $r=1-2^{-1/N}$.
\end{itemize}
\end{corollary}

\medskip
\section{Converting Runs into Run Parts }

The basic idea in the proof of Theorem~2 is to replace a $k$-long run of a given subcomposition $\cvec c$ with a new part
$\kbar$ with $|\kbar|=k|\cvec c|$, and then apply Theorem~1(d) to a new  class $\sA'$ of locally restricted structures with
parts in $\naturals$ as well as {\em run parts} $\kbar$ for $k\ge 1$. Let $\theta$ denote this replacement operation. For example, if $\cvec c=12$ and $\cvec a=12123122121$, then $\theta(\cvec a)={\bar 2}3{\bar 1}2{\bar 1}1$.

\BeginExample[Runs in unrestricted compositions]
Let $\sA$ be the class of all compositions, and $\cvec c$ be a given composition with $\len(\cvec c)=m$ such that it cannot be written in the form $\cvec c=\cvec x\cvec y\cvec x$.  The digraph $D'$ for $\sA'=\theta(\sA)$ is defined as follows.
\begin{itemize}
\item $\sS(D')=\{\eps_s\}$.
\item $\cvec a' \in\sF(D')$ if and only if $\cvec a'\in \SEQ_{<m}(\sP)$, $\cvec a'$ does not contain two consecutive run parts.
\item  $\cvec r'\in \sR(D')$ if and only if $\cvec r'\in \SEQ_{m}(\sP)$, $\cvec r'\ne \cvec c$, $\cvec r'$ does not contain two consecutive run parts.
\item There is an arc from a vertex $\cvec a'$ to a vertex $\cvec b'$ if and only if $\cvec a'\cvec b'$
does not contain $\cvec c$ and does not contain two consecutive run parts.
\end{itemize}
$~~$
\EndExample

Let $\cvec c$ be a given composition which cannot be written in the form $\cvec c=\cvec x\cvec y\cvec x$. It is easy to see that copies of $\cvec c$ in a compositions $\cvec a$ cannot overlap with each other.
Define $\sP=\naturals \cup \{\kbar:k\ge 1\}$ with $|\kbar|=k|\cvec c|$. We have the following
\begin{proposition}\label{bijection}
Let $\sA=\SEQ(\sP;D)$ be a given regular supercritical class of locally restricted compositions with span $m$ such that $\sS(D)=\{\eps_s\}$. Let $\cvec c\in \sR(D)$.  Assume that $D$ contains an arc from $\cvec c$ to itself and  $\cvec c$ cannot be written in the form $\cvec c=\cvec x\cvec y\cvec x$. Let $\sA'=\theta(\sA)$. Then
\begin{itemize}
\item[(a)] $\sA'=\SEQ(\sP;D')$ is a regular supercritical class with $\sP=\naturals \cup \{\kbar:k\ge 1\}$ and some digraph $D'$.
\item[(b)] For each $\cvec a\in {\cal A}$, $|\theta(\cvec a)|=|\cvec a|$, and  $\theta$ is a bijection between ${\cal A}_n$ and ${\cal A'}_n$.
\item[(c)]  For each $\cvec a\in {\cal A}$, Let $R(\cvec a)$ be the maximum run length of $\cvec c$ in $\cvec a$, and $M(\cvec a)$ be the maximum  value of $k$ such that $\kbar$ appears in $\theta(\cvec a)$.  We have
    $R(\cvec a)=M(\cvec a)$.
\end{itemize}
\end{proposition}
\BeginProof (a) We define the digraph $D'$ as follows.
\begin{itemize}
\item $\sS(D')=\{\eps_s\}$.
\item $\cvec a' \in\sF(D')$ if and only if $\cvec a'\in \SEQ_{<m}(\sP)$, $\cvec a'$ does not contain two consecutive run parts and $\theta^{-1}(\cvec a')$ appears at the end of a composition in $\sA$.
\item  $\cvec r'\in \sR(D')$ if and only if $\cvec r'\in \SEQ_{m}(\sP)$, $\cvec r'\ne \cvec c$, $\cvec r'$ does not contain two consecutive run parts,  and $\theta^{-1}(\cvec r')$ appears in a composition in $\sA$.
\item There is an arc from a vertex $\cvec a'$ to a vertex $\cvec b'$ if and only if $\cvec a'\cvec b'$ does not contain two consecutive run parts, does not contain $\cvec c$, and
$\theta^{-1}(\cvec a'\cvec b')$  appears in a composition in $\sA$.
\end{itemize}
 It can be verified that  $\sA'=\SEQ(\sP;D')$ is regular supercritical.  \\
(b) This follows immediately from the definition of $\theta$.\\
(c) Since the copies of $\cvec c$ in $\cvec a$ cannot overlap with each other, $\theta$ induces a bijection between the set of run parts in $\theta(\cvec a)$ and the set of runs in $\cvec a$, which implies $R(\cvec a)=M(\cvec a)$.
\EndProof

\BeginExample[Runs of a subcomposition in $m$-Carlitz compositions]
Let  $\cvec c$ be a sequence of $m$ distinct integers. The digraph $D'$ for $\sA'=\theta(\sA)$ is defined as follows.
\begin{itemize}
\item $\sS(D')=\{\eps_s\}$.
\item $\cvec a' \in\sF(D')$ if and only if $\cvec a'\in \SEQ_{<m}(\sP)$, $\cvec a'$ does not contain two consecutive run parts and $\theta^{-1}(\cvec a')$ is $m$-Carlitz.
\item  $\cvec r'\in \sR(D')$ if and only if $\cvec r'\in \SEQ_{m}(\sP)$, $\cvec r'\ne \cvec c$, $\cvec r'$ does not contain two consecutive run parts, and $\theta^{-1}(\cvec a')$ is $m$-Carlitz.
\item There is an arc from a vertex $\cvec a'$ to a vertex $\cvec b'$ if and only if $\cvec a'\cvec b'$
does not contain $\cvec c$, does not contain two consecutive run parts, and $\theta^{-1}(\cvec a'\cvec b')$ is $m$-Carlitz.
\end{itemize}
$~~$
\EndExample

\BeginExample[Runs of a subcomposition in alternating compositions]
Let $m\ge 2$, $\cvec c=c_1c_2\cdots c_{m}$ be an alternating composition such that $\cvec c^2$ is also alternating. The digraph $D'$ for $\sA'=\theta(\sA)$ is defined as follows.
\begin{itemize}
\item $\sS(D')=\{\eps_s\}$.
\item $\cvec a' \in\sF(D')$ if and only if $\cvec a'\in \SEQ_{<m}(\sP)$, $\cvec a'$ does not contain two consecutive run parts and $\theta^{-1}(\cvec a')$ is alternating.
\item  $\cvec r'\in \sR(D')$ if and only if $\cvec r'\in \SEQ_{m}(\sP)$, $\cvec r'\ne \cvec c$, $\cvec r'$ does not contain two consecutive run parts, and $\theta^{-1}(\cvec a')$ is alternating.
\item There is an arc from a vertex $\cvec a'$ to a vertex $\cvec b'$ if and only if $\cvec a'\cvec b'$
does not contain $\cvec c$, does not contain two consecutive run parts, and $\theta^{-1}(\cvec a'\cvec b')$ is alternating.
\end{itemize}
$~~$
\EndExample

\section{Outline of Proofs}
The proofs are essentially the same as those in \cite{BC09,BCG}. In particular, we will make use of the infinite transfer matrix.
\begin{definition}
[Transfer matrix]
\label{def:transfer}
Let $\SEQ(\sP;D)$ be a class of locally restricted structures as in Definition~\ref{def:LRS}. Let  $\cvec r_1,\cvec r_2,\ldots$ be an ordered list of vertices in $\sR$.
We define the transfer matrix $T(z)$ such that the $(i,j)$th entry of $T(z)$ is $T_{i,j}(z)=z^{|\cvec r_i|+|\cvec r_j|}$ if there is an arc in $D_{\sR}$ from $\cvec r_i$ to $\cvec r_j$; otherwise $T_{i,j}(z)=0$.
\end{definition}

The {\em weight} of an arc $(\cvec v,\cvec w)$ in $D$ is $z^{|\cvec v|+|\cvec w|}$.
We define the weight of a directed walk in $D$ to be the product of the weights of all the arcs in the walk. Let $F_R(z)$ be the ogf for structures in $\SEQ(\sP;D)$ containing at least one recurrent vertex.
It is not difficult to see that $F_R(z^2)$ is the sum of weights of all directed walks in $\sW$ containing at least one recurrent vertex.
We may express $F_R(z)$ in terms of $T(z)$,
the {\em start vector} $\cvec s(z)$ and the {\em finish vector} $\cvec f(z)$, which are defined as follows.
The $i$th component of $\cvec s(z)$ is the sum of weights of all directed walks from $\eps_s$ to $\cvec r_i$, and the $j$th component of
$\cvec f(z)$ is the sum of weights of all directed walks from  $\cvec r_j$ to $\eps_f$.
Since $T^k_{i,j}(z)$ is the sum of weights of all directed walks of length $k$ from $\cvec r_i$ to $\cvec r_j$, we have
$$F_R(z^2)=\cvec s(z)^t\sum_{k\ge 0}T^k(z)\cvec f(z).$$

The following lemma summarizes the results  which are used in the proof of our Theorems 1 and 2. These results are simple extensions of the corresponding results from \cite{BC09} for locally restricted compositions.

\begin{lemma} \label{L1}
Let $\sA=\SEQ(\sP;D)$ be a regular class of locally restricted structures.
\begin{itemize}
\item[(a)]
Let $r$ be the radius of convergence of the generating function $F_R(z)$. Then $r<1$ and it is a simple pole of $F_R(z)$. Moreover  $F_R(z)$ has no other singularity in $|z|\le 1$.
\item[(b)] Let  $p\in \sP$ be a recurrent part and let $X_n$ be the number of occurrences of $ p$ in a random  structure in $\sA_n$.
There are constants $C_i>0$ such that
$$
\Pr(X_n\!<\!C_1n) ~<~ C_2(1+C_3)^{-n}
~~\mbox{for all $n$.}
$$
\item[(c)]
Let $\cvec c$ be a given structure in $\sA$.
There is a constant $B$  such that the
probability that $\cvec c$ occurs in a random structure in $\sA_n$
 is at most $Bnr^{|\cvec c|}$.
\end{itemize}
\end{lemma}

\BeginProofOf{Lemma~\ref{L1}}  Part~(a) is lifted from \cite[Theorem~2]{BC09} and its proof remains exactly the same.\\
 Parts~(b) and (c) are lifted from \cite[Lemma~1]{BC09} and their proofs remain exactly the same.
 \EndProof

\BeginProofOf{Theorem~\ref{thm:main}} Theorem~1(a) follows from Lemma~\ref{L1} as shown in \cite{BC09}. Now the proof of Theorem~1(b--d) is essentially the same as that of \cite[Theorem~1]{BC09}, using  Lemma~1 above, and Lemmas 3, 4 and 5 from \cite{BC09}.
 \EndProof

\smallskip

\begin{lemma}\label{L2}
Let $\sA=\SEQ(\naturals;D)$ be a class of locally restricted compositions (or colored compositions defined in Example~1). Then
$\sA$ is supercritical.
\end{lemma}
\BeginProof Let $G(z)$ be the ogf for all compositions in $\sA$ which don't contain any recurrent vertex.
Since each composition counted by $G(z)$ has at most $K$ parts for some fixed integer $K$,
we have $G(z)\le \sum_{j=0}^{K} P(z)^{j}$ (coefficient-wise). Since the radius of convergence of $P(z)$ is 1, the radius of convergence of
$G(z)$ is at least 1.
\EndProof

\BeginProofOf{Theorem~\ref{thm:Run}} Let $\sP=\naturals\cup \{\kbar:k\ge 1\}$, and $\sA'=\SEQ(\sP;D')$ be the new class defined in Proposition~1.  Since $\theta$ is a bijection between $\sA_n$ and $\sA'_n$, the radius of convergence for $\sA'$ is the same as that for $\sA$. It is easy to see that $\sL=\{\kbar:k\ge 1\}$ is asymptotically free in $\sA'$. Applying Theorem~1(d) with
 $|l_j|=j|\cvec c|$, we obtain
\begin{eqnarray*}
\pr(R_n< k)&=& \prod_{k\le j\le n/|\cvec c|}\pr(\zeta_j(n)=0)\exp\left(-Cnr^{j|\cvec c|}\right) \\
&\sim &  \exp\left(-Cn\frac{r^{k|\cvec c|}}{1-r^{|\cvec c|}}\right).
\end{eqnarray*}
This establishes part~(a).\\
Using $\Ex(R_n)=\sum_{k\ge 1}(1-\pr(R_n\le k-1))$ and the same argument as in \cite[page 25]{BCG} (replacing $r$ with $r^{|\cvec c|}$), we obtain part~(b). We remark that the $+$ sign before $P_0$ should be $-$ in the expression of $f(x)$ in \cite[page 25]{BCG}. The same applies to \cite[Theorem~1(b,c)]{BCG} \\
Part~(c) follows from the same argument as in \cite[page 27]{BCG}, with $r$ being replaced by $r^{|\cvec c|}$.
\EndProof

\BeginProofOf{Corollary~1}
This follows immediately from Theorem~2. The values of $r$ and $C$ are computed as follows. \\
For part~(a), it is clear $r=1/2$. To obtain the value of $C$, we let $F(z)$ be the ogf of compositions containing a marked  $k$-long run of a subcomposition $\cvec c$. We note that $(1-z)/(1-2z)$ is the ogf of all compositions, and so
$z^{|\cvec c|}(1-z)/(1-2z)$ is the ogf of compositions ending with (or starting with) $\cvec c$. Hence we have
$$F(z)=\left(\frac{1-z}{1-2z}-\frac{z^{|\cvec c|}(1-z)}{1-2z}\right)^2z^{k|\cvec c|}=(1-z)^2(1-z^{|\cvec c|})^2z^{k|\cvec c|}(1-2z)^{-2}.$$   It follows from the ``transfer theorem" \cite{FS09} that
$$[z^n]F(z)\sim \frac{n}{4}\left(1-2^{-|\cvec c|}\right)^22^{n-k|\cvec c|}.$$
Hence $C=\frac{1}{2}\left(1-2^{-|\cvec c|}\right)^2$. In particular, $C=1/8$ for the runs of 1, which gives Gafni's result. \\

For part~(b), the value of $r$ can be found, for example, in \cite{BCG}. The expression of $C$ (and $r$) can be derived as follows.  Let $K(z)$ be the ogf of Carlitz compositions, and $F(z)$ be the ogf of Carlitz compositions with a marked $k$-long run of $\cvec c=ab$, where $a\ne b$.  For a given Carlitz composition $\cvec v$, let $K_{\cvec v}(z)$ be the ogf of Carlitz compositions that start with $\cvec v$. Then we have
\begin{eqnarray}
K_a(z)&=&(K(z)-K_a(z))z^a, \label{Ka}\\
K_{ab}(z)&=&z^{a+b}(K(z)-K_b(z)),\nonumber\\
K_{ba}(z)&=&z^{b+a}(K(z)-K_a(z)),\nonumber\\
F(z)&=&z^{k(a+b)}(K(z)-K_b(z)-K_{ab}(z))(K(z)-K_a(z)-K_{ba}(z))\nonumber.
\end{eqnarray}
It follows that
\begin{equation}
F(z)=\frac{1}{(1+z^a)(1+z^b)}z^{k(a+b)}(1-z^{a+b})^2K(z)^2. \label{F}
\end{equation}
The expression of $K(z)$ can be obtained from (\ref{Ka}), which gives
$$K_a(z)=\frac{z^a}{1+z^a}K(z).$$
Summing over $a$, and noting $K(z)=1+\sum_{a\ge 1}K_a(z)$, we obtain
$$K(z)=\frac{1}{1-\sum_{j\ge 1}\frac{z^j}{1+z^j}}.$$
Let $r$ be the smallest positive number satisfying $\sum_{a\ge 1}\frac{z^j}{1+z^j}=1$, it is easy to see that $r$ is a simple pole of $K(z)$ and
$$K(z)=\frac{1}{\sum_{j\ge 1}\frac{jr^j}{(1+r^j)^2}}\frac{1}{1-z/r}+h(z),$$
where $h(z)$ is analytic in $|z|\le r$.
Consequently
$$[z^n]K(z)\sim\frac{1}{\sum_{j\ge 1}\frac{jr^j}{(1+r^j)^2}}r^{-n}.
$$
It follows from (\ref{F}) and the ``transfer theorem" that
\begin{eqnarray*}
~[z^n]F(z)&\sim & \frac{n}{(1+r^a)(1+r^b)}r^{k(a+b)}(1-r^{a+b})^2\left(\frac{1}{\sum_{j\ge 1}\frac{jr^j}{(1+r^j)^2}}\right)^2r^{-n},\\
~\frac{[z^n]F(z)}{[z^n]K(z)}&\sim & \frac{(1-r^{a+b})^2}{(1+r^a)(1+r^b)}\frac{1}{\sum_{j\ge 1}\frac{jr^j}{(1+r^j)^2}}nr^{k(a+b)}.
\end{eqnarray*}
Hence
$$C=\frac{(1-r^{a+b})^2}{(1+r^a)(1+r^b)}\frac{1}{\sum_{j\ge 1}\frac{jr^j}{(1+r^j)^2}}.
$$

The value of $r$ for part~(d) can be found in \cite{BCG}. To obtain the value of $r$ for part~(e), we note that the corresponding part generating function is $P(z)=z(1-z)^{-2}$. Solving the equation $r(1-r)^{-2}=1$, we obtain
$r=\frac{3-\sqrt{5}}{2}$. To obtain the value of $r$ for part~(f), we note that the corresponding part generating function is $P(z)=(1-z)^{-N}-1$. Solving the equation $(1-r)^{-N}=2$, we obtain
$r=1-2^{-1/N}$.
\EndProof

\section{Discussions}

In this paper, we showed how infinite transfer matrix method developed in \cite{BC09} can be applied to enumerate locally restricted regular supercritical sequential structures. Poisson distribution results were derived for sizes of free parts and applications are given for runs of subcompositions in several classes of
compositions. In the upcoming full version of the paper, we plan to extend the results using a more general set up in terms of sized digraphs. Also in Theorem~1, we imposed the condition that the set $\sL$ contains parts of distinct sizes. Such condition is unnecessary for unrestricted sequential structures as shown in \cite{BG14}. It might be possible to relax this condition for some classes of locally restricted sequential structures. For example, when the restriction is size-based, instead of part-based, we may allow the asymptotically free set $\sL$ to contain
several parts of the same size. Finally it might be possible to remove the restriction that the subcomposition $\cvec c$ cannot be written in the form $\cvec c=\cvec x\cvec y\cvec x$ in Theorem~2. The mapping $\theta$ in Proposition~1 is still a bijection provided that the replacement is made at the earliest opportunity. However, the copies of $\cvec c$ in a composition may overlap here, and the maximum run length $R(\cvec a)$ in a composition $\cvec a\in \sA$ may exceed the maximum run part $\bar k$ in $\theta(\cvec a)$. For example, consider the compositions $\cvec c=121$ and $\cvec a=12121121121$. We note that the first two copies of $\cvec c$ overlap at the third position and
$\theta(\cvec a)={\bar 1}21{\bar 2}$. The maximum run part is ${\bar 2}$ here, but the maximum run length of $\cvec c$ in $\cvec a$ is 3.

\vskip 20pt

{}

\end{document}